%% file: main.tex
\label{sec:tracking_prob}

\documentclass[letterpaper, 10 pt, conference]{ieeeconf}  

\input{header.tex}
\SetAlgoSkip{medskip}
\begin{document}
\maketitle
\thispagestyle{empty}
\pagestyle{empty}

\begin{abstract}
Mixed Integer Linear Programs (MILPs) are often used in the path planning of both ground and aerial vehicles.  Such a formulation of the path planning problem requires a linear objective function and constraints, limiting the fidelity of the the tracking of vehicle states.   One such parameter is the state of charge of the battery used to power the vehicle.  Accurate battery state estimation requires nonlinear differential equations to be solved.  This state estimation is important in path planning to ensure flyable paths, however when using MILPs to formulate the path planning problem these nonlinear equations cannot be implemented.  Poor accuracy in battery estimation during the path planning runs the risk of the planned path being feasible by the estimation model but in reality will deplete the battery to a critical level.  To the end of higher accuracy battery estimation within a MILP, we present here a simple linear battery model which predicts the change in state-of-charge (SOC) of a battery given a power draw and duration.  This model accounts for changes in battery voltage due to applied electrical load and changes in battery SOC.  The battery model is presented and then tested against alternate battery models in numerical and in experimental tests.  Further, the effect the proposed linear model has over a simpler SOC estimation on the time-to-solve a resource constrained shortest path problem is also evaluated, implemented in two different algorithms. It is seen that the linear model performs well in battery state estimation while remaining implementable in a Linear Program or MILP, with little affect on the time-to-solve.
\end{abstract}

\section{INTRODUCTION}
Path planning and trajectory planning are problems which show up in many domains, including vehicle routing, Unmanned Aerial Vehicle (UAV) path planning, and general robotics.  These are frequently solved in a discrete manner, where the environment to be traversed is represented as a set of discrete nodes and edges for the agent to travel along.  In applications where a battery is used to power the vehicle, the tracking of battery state-of-charge (SOC) along the path is important to ensure planned paths are feasible on the real system.  If the battery capacity is limited relative to the expected returned path distance, there lies the risk of planning a path which exceeds the vehicle's battery limits but is calculated to be feasible with the battery model used in the planning.  This danger is present in the case of small UAVs (sUAV), hybrid-fuel UAVs which combine battery power with another power source, such as a fuel cell of combustion engine, and electric or hybrid-fuel ground vehicles.  Tracking  battery SOC along a path in a Mixed Integer Linear Program (MILP) requires a linear function for the battery, however battery SOC in reality changes in a nonlinear manner.  To track with a higher fidelity nonlinear battery model, the problem ceases to be linear and thus must be solved with more costly Mixed Integer \textit{Nonlinear} Programming techniques.  However, the cost of increased problem complexity may not be outweighed by the improved battery tracking, depending on the exact problem at hand.  Multi-agent routing, in a MILP form, is often very computationally expensive and quickly becomes infeasible to solve as the problem sizes increase.  If transformed to a nonlinear programming problem, to the end of battery tracking accuracy, the time-to-solve scaling is even harsher.  We present here a simple linear battery model which accounts for voltage changes due to power draw and drop in battery SOC changes.   We compare this to several alternate battery models in both numerical tests and experimental tests in terms of accuracy of state of charge tracking.  Further, this model is implemented in a Resource Constrained Shortest Path Problem (RCSPP) and the time-to-solve is compared to an alternate resource constraint where a constant battery voltage is assumed.   

We are motivated by a prior work in \cite{scott2022power} where noise-aware path planning is done for a hybrid-fuel UAV which combines a generator and battery pack as power sources.  Here a simple equation is used to update the battery SOC where a constant voltage is assumed, so that the change in SOC along an edge is constant regardless of power draw and current SOC.  This is a simplification made to keep the path planning problem easier, where in reality the  change in battery SOC for a given flight leg depends on many factors, two major ones being electrical load and the current battery SOC.  

The primary contributions of this paper are: i) Numerical and experimental testing of a linear battery model, on both single-cell and multi-cell batteries, developed to the end of use in a MILP ii) analysis of the affect the linear model has on the time-to-solve a RCSPP variant as compared to a constant voltage assumption.  The paper is organized as follows.  The battery tracking problem of concern is presented in Section \ref{sec:tracking_prob}.  Section \ref{sec:batt_eq} presents the linear model, along with the alternate models.  Section \ref{sec:num_test} presents the results of a numerical test comparing the linear model with the alternate models on an 18650 battery cell.  Section \ref{sec:exp_test} presents the results of experimental tests on a Lithium Polymer (LiPo) battery.  An evaluation of time-to-solve a RCSPP with and without the linear battery model presented in Section \ref{sec:time_to_solve}.  This RCSPP is solved with two algorithms: i) Branch-and-bound ii) Labeling Algorithm. 

\subsection{Prior Work}\label{sec:priorwork}
To the best of the authors' knowledge, a linear battery model specifically for use in path planning with MILPs has not been presented and tested in the existing literature.  However, there exist numerous models for battery systems.  These are nonlinear and seek to model the battery behavior as accurately as possible.  There is also much work done on real-time estimation of SOC and battery health.  

Battery models are divided into 3 main types: i) Electro-chemical models; ii) Equivalent-circuit models; iii) Data-driven models.  Electro-chemical models seek to model the behavior of the battery at a chemical level, whereas the equivalent-circuit models represent the battery as a circuit system, often a resistor-capacitor (RC) circuit, with behavior that mimics that of the real battery.  Data-driven models develop a numerical model of the battery system based on gathered data.  Battery state estimation techniques either seek to determine the state-of-energy, being the remaining energy that can be pulled from the battery, state-of-charge estimation, power capability, or state-of-health (SOH).  A comprehensive survey of the above estimation and model types is given in \cite{wang2020comprehensive}.  

Tracking battery state in path planning problems has been studied before, primarily in the case of Electric Vehicles (EVs) and sUAVs where the only power source is a battery.  In \cite{di2016coverage}, a model is presented to check if a given path is feasible with regard to battery SOC.  Here, a path is first produced and then checked for battery feasibility after the fact in an online fashion. In \cite{schacht2018path}, the path planning of a UAV is done while accounting for battery SOH.   UAV path planning is considered in \cite{hovenburg2020long} where an algorithm is proposed which includes a sUAV performance model and battery model within the path planning, which was solved using a particle swarm optimization algorithm.  

There is also existing work on Vehicle Routing Problem (VRP) variants applied to electric vehicles (EV).  The Electric VRP (EVRP) is that of finding minimum global cost routes for a fleet of electric vehicles, often with additional constraints such as recharging events.  An extensive survey of the problems and solution approaches is given in \cite{erdelic2019survey}.  There are several studies on hybrid vehicle routing \cite{verma2018electric, doppstadt2016hybrid, vincent2017simulated, hiermann2019routing} where a constant charge/discharge is used for the battery, ignoring the effect battery SOC has on the change in battery SOC for the same edge.  Another study \cite{scott2022power} involves a battery/generator hybrid UAV, where the battery is discharged and charged throughout the flight, where the battery is also tracked using a constant charge/discharge for a given edge.   These problems, all formulated as MILPS, can be improved to have greater accuracy in battery modeling using the model presented in this paper without the significant increase in computational cost when using nonlinear constraints.

\section{Battery Tracking Problem}\label{sec:tracking_prob}
\begin{figure}
    \centering
    \includegraphics[width = 0.48\textwidth]{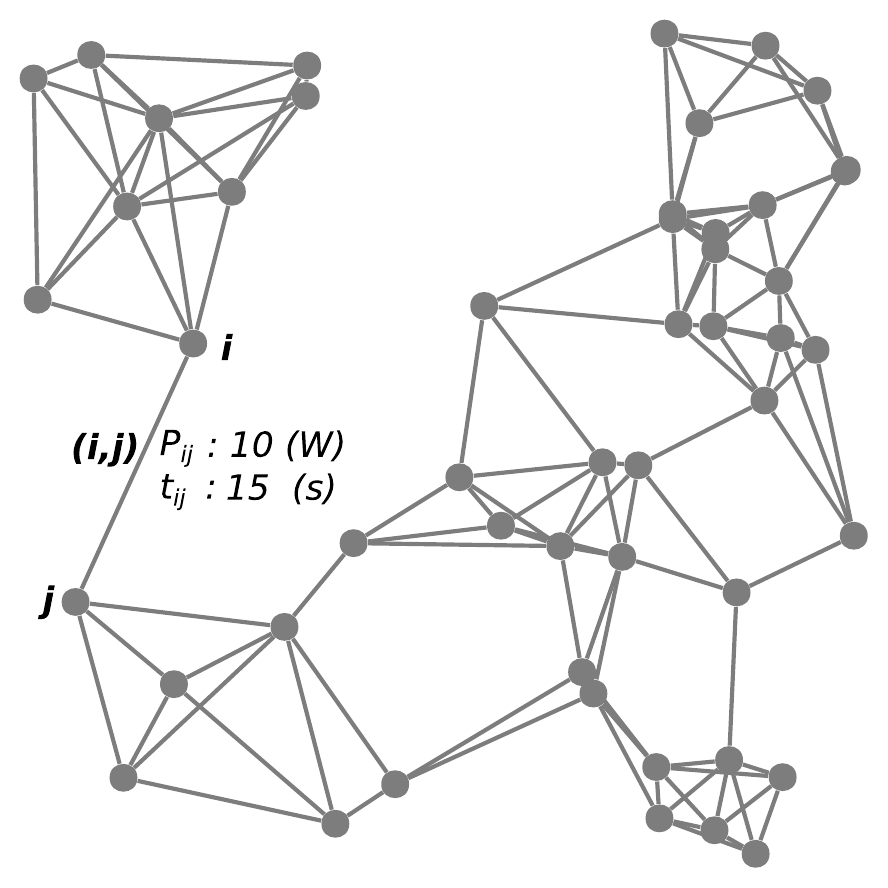}
    \caption{Graph with Power Consumption Parameters}
    \label{figs:param_graph}
\end{figure}
As discussed above, the primary challenge in tracking the battery SOC within a MILP path planning problem arises due to the nonlinear behavior of batteries.  Depending on the chemical principles used for modeling, this nonlinear behavior can be more pronounced.  The path planning, when done on a graph, returns a path, optimal by some metric, as a series of edges.  An example graph is given in Figure \ref{figs:param_graph}, in a case where edges are defined by their power requirement and time to travel along.  Graphs may not be parameterized in this exact manner, but by metrics which track energy or power usage in some alternate manner.  This path planning is often formulated as a MILP, where the battery SOC is tracked by a linear equation that is a function of the graph parameters and decision variables.  An example of such is given in Equation \eqref{battery_update}, a simplified equation from \cite{scott2022power}.

\begin{align}
    &b_j \leq b_i + \Delta(b_i, P_{ij}, t_{ij}) + M(1 - x_{ij}) \quad \forall (i,j) \in E \label{battery_update}  
\end{align}

where $b_i$ is the SOC of the battery at node $i$, $x_{ij}$ is a binary decision variable for usage of the edge from node $i$ to node $j$, and $E$ is the set of all edges in the problem.  The value $M$ is a large constant such that if edge $(i,j)$ is not used in the solution then the constraint is trivially satisfied.  $\Delta(b_i, P_{ij}, t_{ij})$ is the function which determines the change in battery SOC from node $i$ to node $j$ given the power draw $P$ and time duration $t$ for the given edge.  The change in SOC $\Delta$ will be nonlinear.  The aim of the linear model is to estimate the $\Delta$ function while still being implementable in a MILP as a linear constraint.  The above equation is for a shortest path variant, whereas a similar constraint for a VRP or Traveling Salesman Problem variant would require adaptation to the tour format rather than path format.    

\section{Battery Update Equations} \label{sec:batt_eq}
The change in SOC with respect to time is found with Equation \eqref{eq:dsdt}
\begin{align}
    \frac{dS}{dt} = -\frac{I}{C_m} = -\frac{I(S,P)}{C_m} \label{eq:dsdt}
\end{align}
where $S$ is battery SOC normalized by the maximum charge in Coulombs $C_m$ such that at full SOC $S = 1$, $t$ is time, and $I$ is the current through the battery which is positive if the battery is being drained.  Assuming there is a known power draw $P$, the current can be calculated using Equation \eqref{eq:currentpower}.  
\begin{align}
    I(S,P) = \frac{P}{V(S,P)} \label{eq:currentpower}
\end{align}
where $V(S,P)$ is a function of battery SOC and power draw.

Substituting yields the change in SOC with respect to time, given in Equation \eqref{eq:dsdtsubbed}.
\begin{align}
    \frac{dS}{dt} = -\frac{P}{V(S,P) C_m} \label{eq:dsdtsubbed}
\end{align}
Integrating Equation \eqref{eq:dsdtsubbed} gives a function for $S(t)$, given in Equation \eqref{eq:integral}. 
\begin{align}
    S(t) = S_0 - \int_{t_0}^{t} \! \frac{P}{V(S(\tau),P) C_m} \, d\tau \label{eq:integral}
\end{align}
where a positive power value $P$ is that which drains the battery.  This can be solved numerically with methods such as Runge-Kutta techniques or similar.  The solution obtained from the numerical methods can be compared to the linear model estimation, which uses a single time-step in calculation of $S(t)$.  For comparison against the linear model, three alternate battery models are used.  The first is a simple Ohmic drop model and the second is a 1st order RC model, both from \cite{hu2012comparative}.  These are models (2) and (7) respectively as referred to in the cited study.  The Ohmic drop model was chosen due to its simplicity, comparable to that of the linear model presented here.  The 1st order RC model was chosen as its accuracy was comparable to more complex models presented in \cite{hu2012comparative}, where the higher complexity models saw diminishing returns in accuracy improvements.  The third alternate model used is one which assumes constant voltage.  Here, this is taken to be the nominal voltage of the battery types testing in Section \ref{sec:num_test} and \ref{sec:exp_test}.

\subsection{Simple Ohmic Drop Model} 
The open-circuit voltage (OCV) of a battery can be obtained experimentally and utilized by a look-up table to find OCV for a given SOC.  However, when a load is applied to the battery, the actual voltage supplied by the battery drops.  This effect, called Ohmic drop, increases with increasing load on the battery.  For higher loads, the battery is supplying a lower voltage and thus, for the same power draw, a larger current.  Therefore, the change in SOC depends on both the current SOC and the power.  The simple Ohmic drop model is defined as:
\begin{align}
    & V(S,I) = OCV(S) - I R_0 \label{eq:simpleohmic} \\
    & P = V(S,I) I \label{power_k}
\end{align}
where $R_0$ is the internal resistance of the battery. The current $I$ is assumed to be a function of power and voltage, as in our case there is a constant power draw through which the current $I$ is defined, whereas in most cases the current is assumed to be known directly.  Substituting \eqref{power_k} into \eqref{eq:simpleohmic} gives the following where voltage is now a function of SOC and power draw:
\begin{align}
    & V(S,P) = \frac{OCV(S) + \sqrt{OCV(S)^2 - 4 P R_0  } }{2} \label{eq:simpleohmicsolved}
\end{align}
Equation \eqref{eq:simpleohmicsolved} is the simple Ohmic drop model when constant power is used. Here, $V(S,P)$ will always have a real, positive solution so long as the current is not so large to cause a negative voltage. Thus, the cases of complex solutions are ignored.  

\subsection{1st Order RC Model}      
The first-order model is defined in Equations \eqref{eq:RC1} and \eqref{eq:RC2}.
\begin{align}
    & V_k = OCV(SOC_k) - I_k R_0 - U_{k} \label{eq:RC1}\\
    & U_{k+1} = \exp(-\tfrac{\Delta t}{\tau}) U_{k} + R_1 [1 -  \exp(-\tfrac{\Delta t}{\tau})   ] \frac{P_k}{V_k} \label{eq:RC2}
\end{align}
where $\tau$ is the time constant of the RC circuit which models the battery behavior, $R_0$ is the internal resistance of the battery, $k$ is the numerical step, $R_1$ and $\tau$ are the resistance and time constant of the simulated circuit, and $U_k$ is a secondary variable tracking the hysteresis effects over time.

\subsection{Nominal Voltage Only}
The fourth estimation model used is one which uses the nominal battery voltage exclusively, and determines change in SOC based off this nominal voltage, power draw, and time duration of the load.  This type of model is used often in path planning, as discussed in Section \ref{sec:priorwork}.  This estimation model is solved using a single time step as a fair comparison to the linear model, as both can be implemented in a MILP in this manner. The update equation is given in Equation \eqref{eq:nominalmodel}, where $V_{nom}$ is the nominal battery voltage.

\begin{align}
 &S(t) =  S_0 - \frac{P (t - t_0)}{C_m V_{nom} }  \label{eq:nominalmodel}
\end{align}

\subsection{Linear Model}\label{sec:linear_model} 
The linear model to implement Equation \eqref{eq:integral} into a MILP constraint is now presented.  A linear fit is made utilizing the simple Ohmic drop in Equation \eqref{eq:simpleohmic}.  However, since the change in SOC depends on $\frac{1}{V}$, a linear fit will be made of the inverse of the voltage given in Equation \eqref{eq:simpleohmic}. The linear model which accounts for the power draw and current SOC is defined in Equation \eqref{eq:linmodel}.
\begin{align}
 &S(t) =  S_0 -   P (A S_0 + B P + C) (t - t_0) / C_m \label{eq:linmodel}
\end{align}
where (A,B,C) are the fit parameters of the linear model. This is a linear  fit of $\frac{1}{V(S, P)}$ where $V(S,P)$ is calculated with Equation \eqref{eq:simpleohmicsolved}.   Note that over the time period $(t_0,t)$, the voltage will change as the battery is discharged. This is addressed in the RC model and simple Ohmic drop model via the Runge-Kutta techniques to solve Equation \eqref{eq:integral}.  As the linear model is to be implemented in a MILP, such techniques cannot be used and instead a single time-step is used, where the voltage, which changes over the time period, is taken to be constant as the voltage at time $t_0$.  

\begin{figure}
    \centering
    \includegraphics[width = 0.48\textwidth]{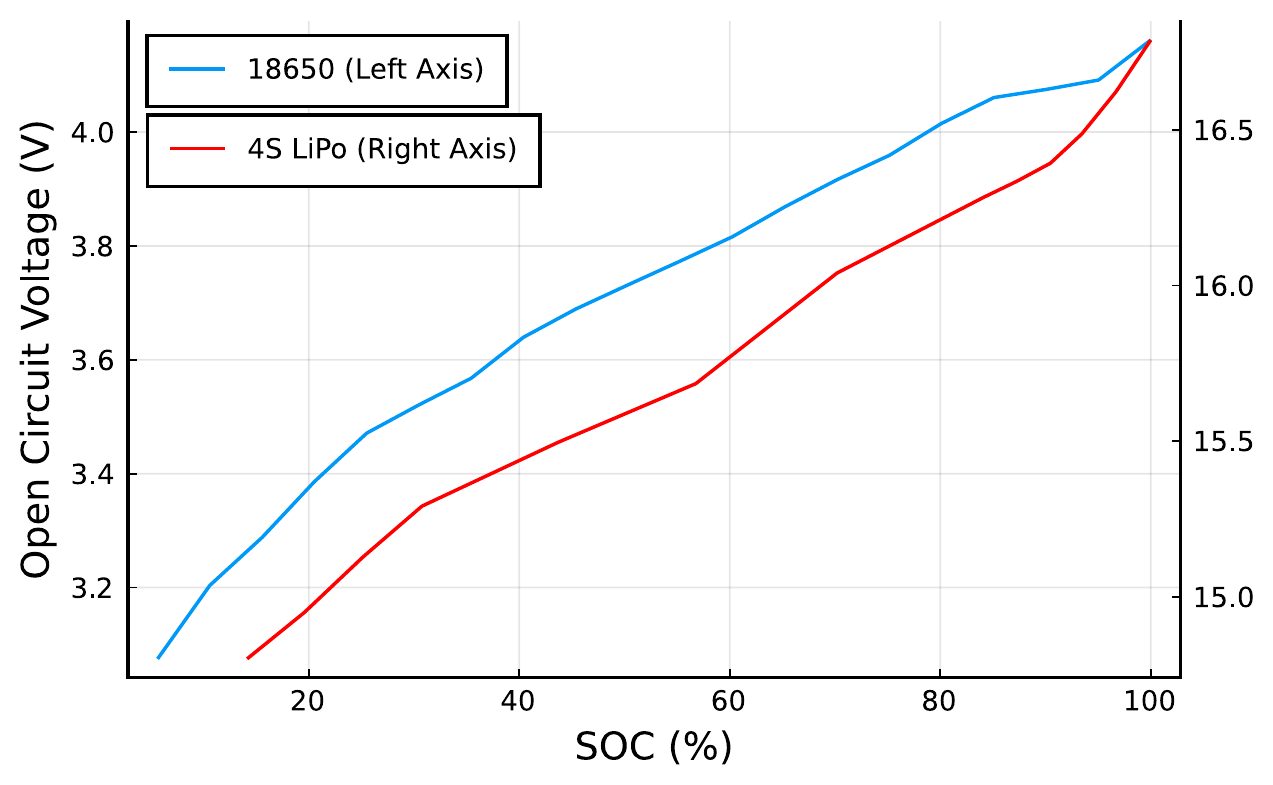}
    \caption{OCV Curve - 4S LiPo and 18650 Cell}
    \label{fig:OCV}
\end{figure}
\begin{figure}
    \centering
    \includegraphics[width = 0.48\textwidth]{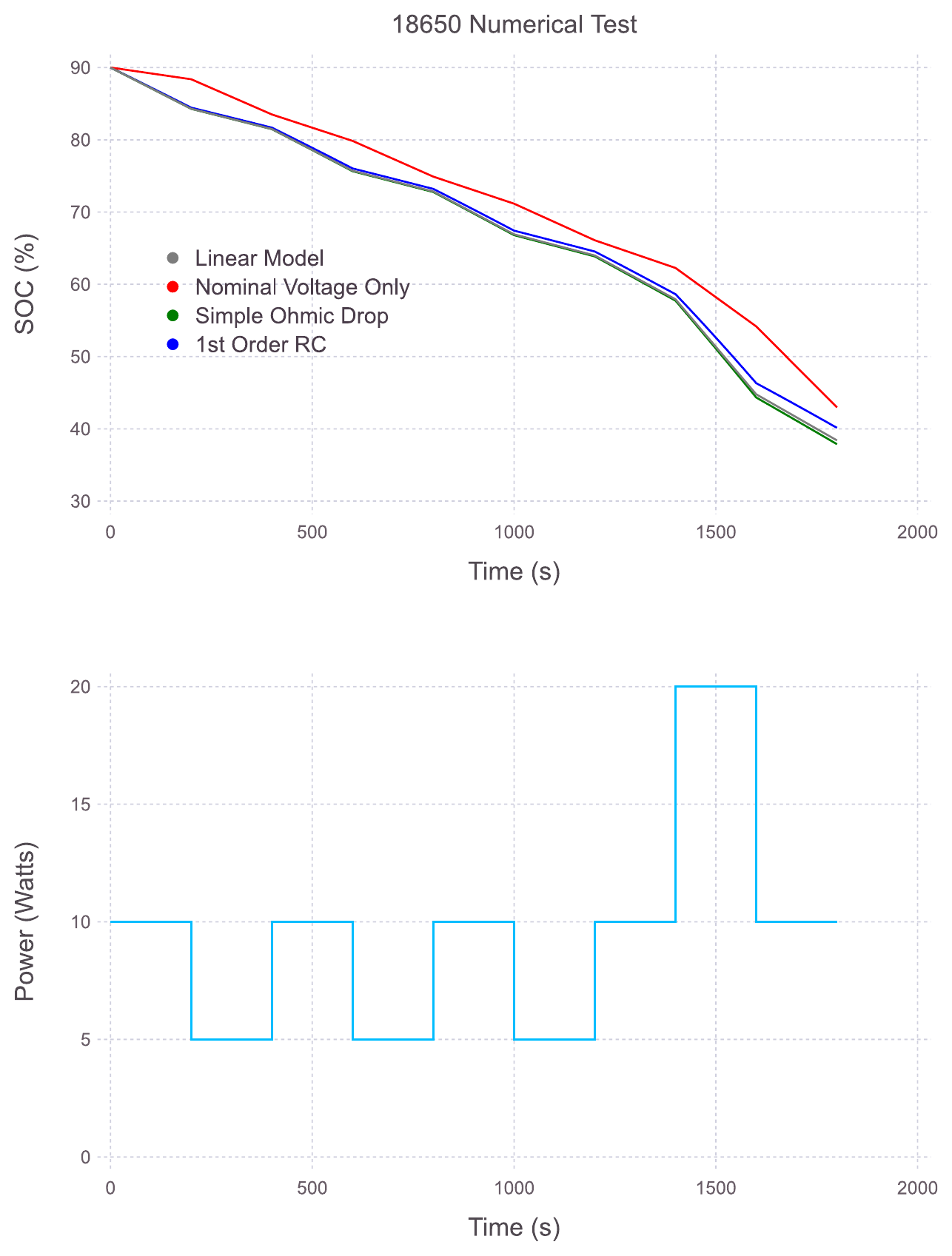}
    \caption{Numerical Test - 18650 cell}
    \label{fig:18650_pulse}
\end{figure}
\begin{figure}
    \centering
    \includegraphics[width = 0.48\textwidth]{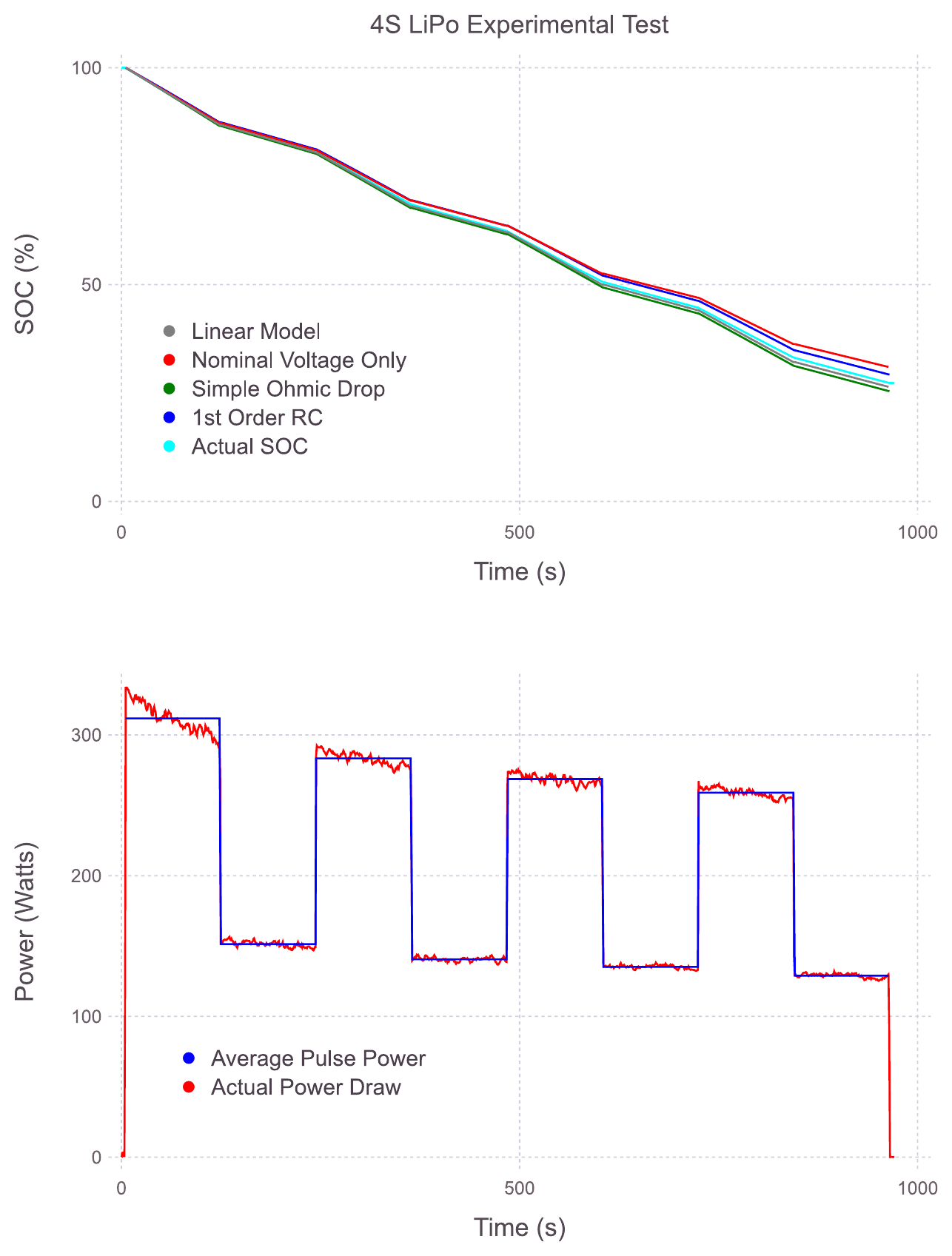}
    \caption{Experimental Test - 4S LiPo}
    \label{fig:LiPo_pulse}
\end{figure}

\section{Numerical Tests} \label{sec:num_test}
Here we present the results of a numerical simulation for a single 18650 battery cell, where the constant-power flight legs of the path-planning problem are simulated.  In the path-planning problem, the battery SOC is to be tracked along a path, where constant power is drawn for some known duration.  The path, in this form, is a queue of these flight legs or edges, and thus can be represented as a series of constant power draws from the battery for some duration.  An 18650 lithium-ion (Li-Ion) battery with 2500 mAh capacity is tested, with the OCV-SOC look up table taken from \cite{elmahdi2021fitting}.  Plots of the OCV-discharge curves are given in Figure \ref{fig:OCV}.  In each set of tests the 4 estimation models from Section \ref{sec:batt_eq} are used: i) simple Ohmic drop model in Equation \eqref{eq:simpleohmic} ii) 1st Order RC model iii) Linear Model presented iv) Nominal Voltage only   

The results of the simulation for the 18650 Li-Ion cell are given in Figure \ref{fig:18650_pulse}.  Here, the power pulses are applied such that the battery is nearly depleted by the end of the simulation, so that the models are tested over nearly the entire SOC of the battery.  It can be seen that all models perform similarly, with a few percent-SOC difference in the final estimations.  Because the power applied is constant, and thus the current changing very slowly only as a result of change in battery voltage as the battery is depleted, the transient effects of the battery are minimal in this application.  Cases were the power changes more frequently may show improved performance of the 1st order RC model relative to the others as the transient effects will have a greater influence. 

\section{Experimental Tests} \label{sec:exp_test}
The results of using the models to predict battery SOC as compared to real data from a 4S LiPo battery are presented in this section.  A 4S LiPo was used, with 5000 mAh capacity, and a 50C maximum discharge rate.  Measurements for power, voltage, and current were used with an 1585 RC Benchmark test stand.  The battery load was applied via a motor and propeller, with the motor controlled by an Electronic Stability Control (ESC) Pulse-Width Modulation (PWM) setting.  In these tests, a constant ESC frequency was applied to achieve approximately the desired power draw.  As seen in Figure \ref{fig:LiPo_pulse}, the power draw was not constant, including both the noise and slow reduction in power draw over time.  The drop in power is due to the constant ESC setting rather than power draw.  As the battery discharged and heated up throughout the test, the voltage decreases, which then reduces power draw for the same ESC PWM setting.  From the battery data, a constant power draw was calculated for each "pulse".  This constant power draw, along with time duration of that specific pulse, was given to each model to predict the battery discharge along that pulse.

The results of the model predictions and the actual battery SOC are given in Figure \ref{fig:LiPo_pulse}.  It can be seen that the models compare similarly among themselves, as was the case in the simulation of the 18650 battery in Section \ref{sec:num_test}. Between the linear model and the actual battery SOC there was a 0.84\% SOC difference.

The linear model, due to the limitations of being implemented in a MILP, only uses a single time step to calculate change in state of charge.  When compared to the prediction model which also uses a single time step but only uses the nominal voltage of the battery to predict change in state of charge, the linear model is more accurate while remaining implementable in a MILP for path planning or mission planning.  The other models tested here, described in Section \ref{sec:batt_eq}, used Runge-Kutta with 100 time-steps to solve Equation \eqref{eq:integral} whereas the linear model uses a single time step.  However, despite the drastic difference in time-steps, the linear model is still able to predict battery SOC  with high accuracy for the path planning battery prediction problem in both the numerical test and on the experimental data.

\section{MILP Path Planning Test}\label{sec:time_to_solve}
By design, the linear model is implementable in a MILP for path planning, however increased complexity in MILP constraints can have a significant effect on computation time.  This effect can also vary between algorithms used to solve the problem.  We implement the linear battery model in a RCSPP \cite{joksch1966shortest} and compare the time to solve with the same problems using nominal battery voltage only.  The MILP formulation is given in Equations \eqref{eq:cost}-\eqref{eq:RCSPP_naive}.
\begin{align}
    & \min_{x} \sum_{i} \sum_{j} ( C_{ij} x_{ij}) \label{eq:cost}\\ 
    & \sum_{j} x_{Sj} = 1 \label{eq:start}\\
    & \sum_{j} x_{jF} = 1 \label{eq:end}\\
    & \sum_{j}x_{ij} - \sum_{j}x_{ji} = 0  && \forall i \in N \setminus \{S,F\} \label{eq:cont} \\
    & 0 \leq b_j \leq B_m && \forall j \in N \setminus {S} \label{eq:batt1}\\ 
    & b_S = B_0 \label{eq:batt2}
\end{align}
where $x_{ij}$ is a binary decision variable which is $1$ if edge from node $i$to node $j$ is used in the solution and 0 otherwise, $b_i$ is the battery SOC  at node $i$, and $S$ and $F$ are the start and goal nodes respectively.  The constraints to track battery state are given as follows:
\begin{align}
    & b_j \leq b_i - P_{ij} t_{ij} (A \hspace{1pt} b_i + B \hspace{1pt}P_{ij}  + C  ) \nonumber \\ 
    & \qquad \quad + M(1 - x_{ij}) \quad && \forall (i,j) \in E \label{eq:RCSPP_lin} \\
    & b_j \leq b_i - \frac{P_{ij} t_{ij}}{V_{nom}}  + M(1 - x_{ij}) \quad && \forall (i,j) \in E \label{eq:RCSPP_naive}
\end{align}
where $P_{ij}$ is the power draw, $t_{ij}$ is the time to travel along the edge, $(A,B,C)$ are the linear model parameters, and $V_{nom}$ is the nominal voltage of the battery.  
 The linear battery model is used as the resource constraint as implemented with Equation \eqref{eq:RCSPP_lin}.  The alternate resource constraint, using the nominal voltage model, is given in \eqref{eq:RCSPP_naive}.

\begin{figure}
    \centering
    \includegraphics[width = 0.47\textwidth]{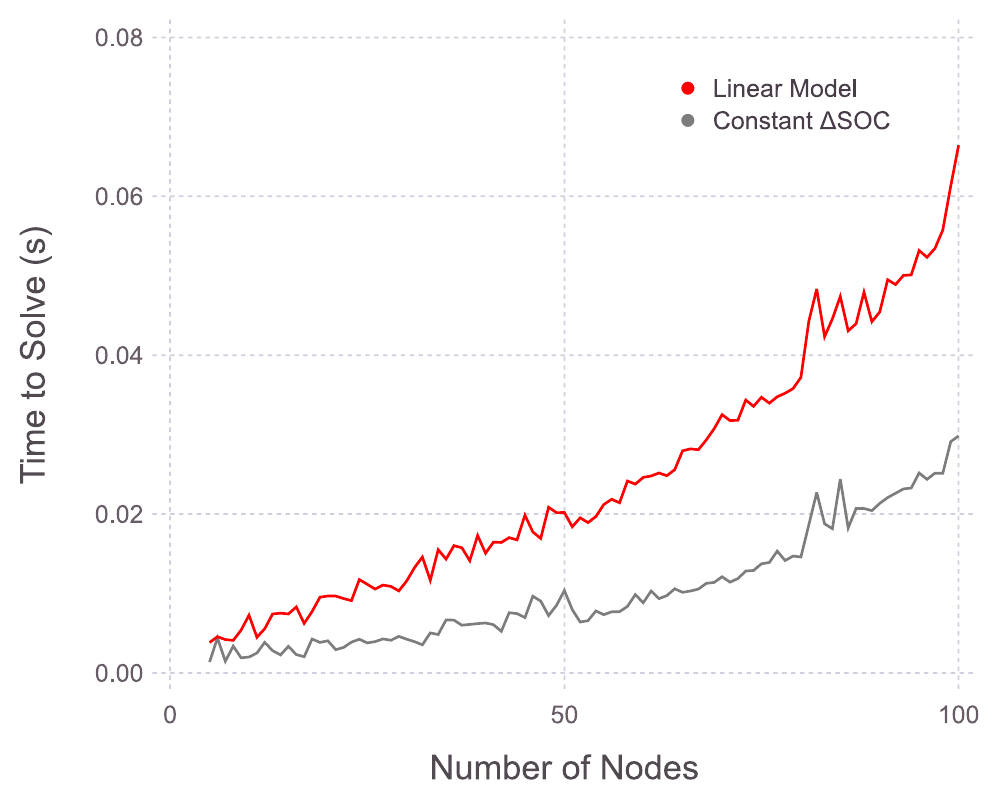}
    \caption{Time-to-Solve with Branch-and-Bound - Battery Model Comparison}
    \label{fig:timetosolve_CPLEX}
\end{figure}

\begin{figure}
    \centering
    \includegraphics[width = 0.47\textwidth]{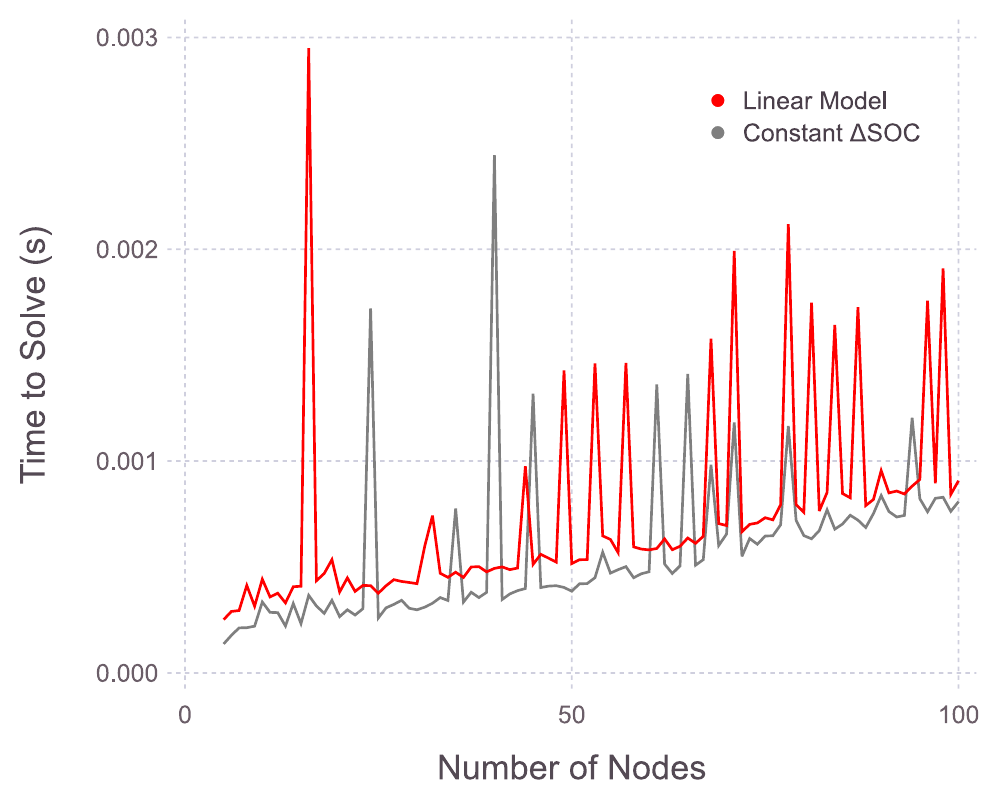}
    \caption{Time-to-Solve with Labeling Algorithm- Battery Model Comparison}
    \label{fig:timetosolve_label}
\end{figure}

The effect on running time when using the linear model over the nominal voltage only model is evaluated when solving with 2 different algorithms: i) Branch-and-bound ii) Labeling algorithm.  The MILP solver in IBM CPLEX v12.9 is used to solve with brand-and-bound and the labeling algorithm presented in \cite{scott2022power} is modified for this RCSPP problem.  Problems are made from 5 nodes up to 100 nodes, with 30 instances generated for each problem size. Nodes are randomly placed in a 2D space and each node is connected to the 4 closest nodes, by euclidean distance.  Edge costs are set as the Euclidean distance between nodes.  The battery parameters for the problem are defined to simulate the 4S LiPo battery described in Section \ref{sec:exp_test}.  The averaged results for time-to-solve with the branch-and-bound is given in Figure \ref{fig:timetosolve_CPLEX} and for the labeling algorithm in Figure \ref{fig:timetosolve_label}.  

The time to solve is increased in both the branch-and-bound and the labeling algorithm when using the linear model.  However, this increase in computation time is more pronounced in the branch-and-bound as compared to the labeling algorithm.  This is expected, as the constraint complexity can have a significant effect on time-to-solve when using a branch-and-bound method.  The labeling algorithm, which utilizes dynamic programming, has the computation time primarily affected by the number of feasible, undominated paths and sub-paths that exist within the problem.  While the linear model increases the battery constraint complexity, it has little effect on the number of feasible, undominated paths within the problem.  Therefore, the battery constraint is expected to have little effect on the computational time when using the labeling algorithm.  The spikes seen in the time-to-solve with the labeling algorithm are expected, where infeasible and loosely constrained problems require many more iterations to complete compared to the average number of iterations, for the same problem size.  A presentation of labeling algorithms and their use in constrained shortest path problem variants is given in \cite{desrosiers1995time}.

\section{Conclusions and Future Work}
The linear battery model was presented to the end of accurate battery SOC estimation for MILP path planning problems concerned with EVs and battery-powered UAVs. The model was applied to an 18650 battery cell and a 4S LiPo battery pack where a current was applied to simulate the series of steady-state flights along graph edges in the path planning problem.  Beyond the accuracy testing, changes to the time-to-solve a RCSPP were also studied.  It was seen that the linear model provides greater accuracy in SOC estimation over assuming a constant nominal voltage with little effect on computational time, while remaining a linear constraint able to be applied in a MILP.  Future work includes further testing on the linear model under different circumstances, and integration with a power or energy consumption model to accurately predict power usage along an edge. An accurate power model is important for path planning for the same reason as the battery model, that being the need to avoid paths which are feasible by the battery and power models but in actuality are not feasible as the battery will drain to a critical level before reaching the goal.  
\bibliographystyle{IEEEtran}
\bibliography{bibfile}
\end{document}

%% file: header.tex
\IEEEoverridecommandlockouts                              

\usepackage{graphics} 
\usepackage{graphicx}
\usepackage{epsfig} 
\usepackage{times} 
\usepackage{amsmath,amssymb} 

\usepackage[linesnumbered,ruled]{algorithm2e}
\usepackage[justification=centering]{caption}
\usepackage{flushend}
\usepackage{subfloat}

\allowdisplaybreaks

\title{\LARGE \bf
Development of Linear Battery Model for Path Planning with Mixed Integer Linear Programming: Simulated and Experimental Validation
}

\author{Drew Scott$^{1}$, Satyanarayana G. Manyam$^{2}$, David W. Casbeer$^{3}$, Manish Kumar$^{4}$,
\\ Isaac E. Weintraub$^{5}$, Michael J. Rothenberger$^{6}$ 
\thanks{$^{1}$Department of Mechanical and Materials Engineering, University of Cincinnati 
        {\tt\small scott2dd@mail.uc.edu}}%
\thanks{$^{2}$Research Scientist, Infoscitex corp., a DCS company, Dayton OH, 45431
        {\tt\small msngupta@gmail.com}}%
\thanks{$^{3}$Technical Area Lead, Cooperative \& Intelligent Control, Control Science
Center, Air Force Research Laboratory, WPAFB, OH, 45433
        {\tt\small david.casbeer@us.af.mil}}%
\thanks{$^{4}$ \noindent Professor in Department of Mechanical and Materials Engineering, University of Cincinnati {\tt\small kumarmu@ucmail.uc.edu}}
\thanks{$^{5}$ \noindent Electronics Engineer, Control Science Center, Air Force Research Laboratory, WPAFB, OH, 45433 {\tt\small isaac.weintraub.1@us.af.mil}}%
\thanks{$^{6}$ \noindent Mechanical Engineer, Flight Systems Integration Branch, Air Force Research Laboratory, WPAFB, OH, 45433}%
\thanks{Distribution Statement A: Approved for Public Release; Distribution is Unlimited. PA\# AFRL-2022-4598}
\thanks{This  work  has  been  supported  in  part  by  AFOSR  LRIR  No.  21RQ-COR084.}
\thanks{We acknowledge the support from Bryan Brown and Austin Wessels at the UAV MASTER Lab at University of Cincinnati}
\thanks{This material is based on research sponsored by the Air Force Research Laboratory and Strategic Council for Higher Education under agreement FA8650-19-2-9300.  The U.S. Government is authorized to reproduce and distribute reprints for Governmental purposes notwithstanding any copyright notation thereon. The views and conclusions contained herein are those of the authors and should not be interpreted as necessarily representing the official policies or endorsements, either expressed or implied by the Strategic Council for Higher Education and the Air Force Research Laboratory (AFRL) or the U.S. Government.}
}

%% file: main.bbl
\begin{thebibliography}{10}
\providecommand{\url}[1]{#1}
\csname url@samestyle\endcsname
\providecommand{\newblock}{\relax}
\providecommand{\bibinfo}[2]{#2}
\providecommand{\BIBentrySTDinterwordspacing}{\spaceskip=0pt\relax}
\providecommand{\BIBentryALTinterwordstretchfactor}{4}
\providecommand{\BIBentryALTinterwordspacing}{\spaceskip=\fontdimen2\font plus
\BIBentryALTinterwordstretchfactor\fontdimen3\font minus
  \fontdimen4\font\relax}
\providecommand{\BIBforeignlanguage}[2]{{%
\expandafter\ifx\csname l@#1\endcsname\relax
\typeout{** WARNING: IEEEtran.bst: No hyphenation pattern has been}%
\typeout{** loaded for the language `#1'. Using the pattern for}%
\typeout{** the default language instead.}%
\else
\language=\csname l@#1\endcsname
\fi
#2}}
\providecommand{\BIBdecl}{\relax}
\BIBdecl

\bibitem{scott2022power}
D.~Scott, S.~G. Manyam, D.~W. Casbeer, M.~Kumar, M.~J. Rothenberger, and I.~E.
  Weintraub, ``Power management for noise aware path planning of hybrid
  {UAVs},'' in \emph{American Control Conference (ACC)}, Atlanta, GA, USA,
  2022, pp. 4280--4285.

\bibitem{wang2020comprehensive}
Y.~Wang, J.~Tian, Z.~Sun, L.~Wang, R.~Xu, M.~Li, and Z.~Chen, ``A comprehensive
  review of battery modeling and state estimation approaches for advanced
  battery management systems,'' \emph{Renewable and Sustainable Energy
  Reviews}, vol. 131, p. 110015, 2020.

\bibitem{di2016coverage}
C.~Di~Franco and G.~Buttazzo, ``Coverage path planning for {UAVs}
  photogrammetry with energy and resolution constraints,'' \emph{Journal of
  Intelligent \& Robotic Systems}, vol.~83, no.~3, pp. 445--462, 2016.

\bibitem{schacht2018path}
R.~Schacht-Rodr{\'\i}guez, J.-C. Ponsart, C.-D. Garc{\'\i}a-Beltr{\'a}n, C.-M.
  Astorga-Zaragoza, D.~Theilliol, and Y.~Zhang, ``Path planning generation
  algorithm for a class of {UAV} multirotor based on state of health of lithium
  polymer battery,'' \emph{Journal of Intelligent \& Robotic Systems}, vol.~91,
  no.~1, pp. 115--131, 2018.

\bibitem{hovenburg2020long}
A.~R. Hovenburg, F.~A. de~Alcantara~Andrade, R.~Hann, C.~D. Rodin, T.~A.
  Johansen, and R.~Storvold, ``Long-range path planning using an aircraft
  performance model for battery-powered {sUAS} equipped with icing protection
  system,'' \emph{IEEE Journal on Miniaturization for Air and Space Systems},
  vol.~1, no.~2, pp. 76--89, 2020.

\bibitem{erdelic2019survey}
T.~Erdeli{\'c} and T.~Cari{\'c}, ``A survey on the electric vehicle routing
  problem: variants and solution approaches,'' \emph{Journal of Advanced
  Transportation}, vol. 2019, pp. 1--48, 05 2019.

\bibitem{verma2018electric}
A.~Verma, ``Electric vehicle routing problem with time windows, recharging
  stations and battery swapping stations,'' \emph{EURO Journal on
  Transportation and Logistics}, vol.~7, no.~4, pp. 415--451, 2018.

\bibitem{doppstadt2016hybrid}
C.~Doppstadt, A.~Koberstein, and D.~Vigo, ``The hybrid electric
  vehicle--traveling salesman problem,'' \emph{European Journal of Operational
  Research}, vol. 253, no.~3, pp. 825--842, 2016.

\bibitem{vincent2017simulated}
F.~Y. Vincent, A.~P. Redi, Y.~A. Hidayat, and O.~J. Wibowo, ``A simulated
  annealing heuristic for the hybrid vehicle routing problem,'' \emph{Applied
  Soft Computing}, vol.~53, pp. 119--132, 2017.

\bibitem{hiermann2019routing}
G.~Hiermann, R.~F. Hartl, J.~Puchinger, and T.~Vidal, ``Routing a mix of
  conventional, plug-in hybrid, and electric vehicles,'' \emph{European Journal
  of Operational Research}, vol. 272, no.~1, pp. 235--248, 2019.

\bibitem{hu2012comparative}
X.~Hu, S.~Li, and H.~Peng, ``A comparative study of equivalent circuit models
  for {Li-ion} batteries,'' \emph{Journal of Power Sources}, vol. 198, pp.
  359--367, 2012.

\bibitem{elmahdi2021fitting}
F.~Elmahdi, L.~Ismail, and M.~Noureddine, ``Fitting the {OCV-SOC} relationship
  of a battery lithium-ion using genetic algorithm method,'' in \emph{E3S Web
  of Conferences}, vol. 234.\hskip 1em plus 0.5em minus 0.4em\relax EDP
  Sciences, 2021, p. 00097.

\bibitem{joksch1966shortest}
H.~C. Joksch, ``The shortest route problem with constraints,'' \emph{Journal of
  Mathematical analysis and applications}, vol.~14, no.~2, pp. 191--197, 1966.

\bibitem{desrosiers1995time}
J.~Desrosiers, Y.~Dumas, M.~M. Solomon, and F.~Soumis, ``Time constrained
  routing and scheduling,'' \emph{Handbooks in operations research and
  management science}, vol.~8, pp. 35--139, 1995.

\end{thebibliography}
